\documentclass{article}

\usepackage{amsmath,amssymb}
\newcommand{\R}{\mathbb{R}}
\newcommand{\N}{\mathbb{N}}
\newcommand{\Z}{\mathbb{Z}}
\newcommand{\y}{{\boldsymbol{y}}}
\newcommand{\Y}{{\boldsymbol{Y}}}
\newcommand{\p}{\hat{P}}
\newtheorem{theorem}{Theorem}
\newtheorem{lemma}[theorem]{Lemma}
\newtheorem{proposition}[theorem]{Proposition}

\DeclareMathOperator{\var}{Var}

\begin{document}
\title{On the Spectral Gap for Convex Domains   }
\date{July 6, 2006}
\author{Burgess Davis\\
Department of Mathematics,\\
Purdue University,\\
150 N. University Street,\\
West Lafayette, IN 47907--2067 \\
E-mail: bdavis@stat.purdue.edu\and Majid Hosseini\\
Department of Mathematics,\\
State University of New York at New Paltz,\\
1 Hawk Drive. Suite 9,\\
New Paltz, NY 12561--2443\\
E-mail: hosseinm@newpaltz.edu} \maketitle
\begin{abstract} Let $D$ be a convex planar domain, symmetric
about both the $x$- and $y$-axes, which is strictly contained in
$(-a,a)\times (-b,b)=\Gamma$. It is proved that, unless $D$ is a
certain kind of rectangle, the difference (gap) between the first
two eigenvalues of the Dirichlet Laplacian in $D$ is strictly
larger than the gap for $\Gamma$.  We show how to give explicit
lower bounds for the difference of the gaps.
\end{abstract}

\section{Introduction}
Let $0<\lambda_1^{\Omega}<\lambda_2^{\Omega}$ be the first two
eigenvalues of the Dirichlet Laplacian for the bounded planar
domain $\Omega$. This paper is concerned with the spectral gap
$\lambda_2^{\Omega}-\lambda_1^{\Omega}$ of $\Omega$. The gap is
the rate at which the Dirichlet heat kernel
$p_t^{\Omega}(x,\cdot)$, normalized to have integral one,
converges to the first eigenfunction, also normalized to have
integral one, where convergence here can mean $L^1$~convergence,
$L^2$~convergence, or pointwise convergence. We note the paper
\cite{swyy} says the gap ``needs no motivation.'' Other papers
concerned with gaps of convex planar domains include
\cite{bankro,banmen,vdb,davis,Ling93,swyy,smits,Wang00,YuZhong86}.

The gap of $(-a,a)\times (-b,b)=\Gamma$ is $3\pi^2/4\max(a,b)^2$.
Davis proved in \cite{davis} that if $D$ is doubly symmetric and
convex, and contained in $\Gamma$, the gap of $D$ is no smaller
than the gap of $\Gamma$. Neither the proof of this in
\cite{davis} nor subsequent proofs (see \cite{banmen},
 \cite{bankro}, and \cite{draghici}) give results about strict inequality. Now
$(-c,c)\times (-b,b)$ is strictly contained in $\Gamma$ if $c<a$,
but if $b\geq a$ it has the same gap, $3\pi^2/4b^2$, as $\Gamma$.
We show such rectangles are the only exceptions.

\begin{theorem}\label{th:1}
If $D$ is convex and symmetric about both the $x$-  and $y$-axes
and strictly contained in $(-a,a)\times(-b,b)$, and is not a
rectangle of the form $(-c,c)\times (-b,b)$ or $(-a,a)\times
(-c,c)$, then the gap of $D$ exceeds the gap of $(-a,a)\times
(-b,b)$.
\end{theorem}

 In common with previous work in
\cite{banmen}, \cite{draghici}, and \cite{davis}, our proof of
Theorem~\ref{th:1} uses the consequence of a theorem of Payne
\cite{payne} that for $D$ as in Theorem~\ref{th:1}, either the
intersection of the $x$-axis with $D$ or the intersection of the
$y$-axis with $D$ is a nodal line for a second eigenfunction. Thus
a second eigenfunction  of $D$ is the first eigenfunction of
either the right or the top half of $D$. The first eigenvalue is
the rate of decay of the heat kernel, and so the proof given in
the next section, that the following proposition implies
Theorem~\ref{th:1}, is quick.  For any set $A\subset\R^2$, let
$A^+=\{(x,y)\in A \mid x>0\}$.

\begin{proposition}\label{pr:2}
Let $D$ be a  convex domain and symmetric about both the $x$- and
the $y$-axes such that $(a,0)$, and $(0,b)$ are boundary points of
$D$. Suppose $\Gamma=(-a,a)\times (-b,b)$ strictly contains $D$.
Then if $z_0=(x_0, y_0)\in D^+$,
\begin{equation}\label{eq:1.1}
\frac{\int_{D^+} p_t^{D^+} (z_0,z)\,dz}{\int_D p_t^D(z_0,z)\,dz} =
o(1) \frac{\int_{{\Gamma}^+} p_t^{{\Gamma}^+}
(z_0,z)\,dz}{\int_{\Gamma} p_t^{\Gamma}(z_0,z)\,dz} \quad \mbox{
as $t\rightarrow \infty$.}
\end{equation}
\end{proposition}

 It is not hard to show that
Proposition~\ref{pr:2} implies  that $o(1)$  in fact decreases
exponentially as $t\rightarrow \infty$.  Our proof yields
estimates on $o(1)$ which depend on the shape of $D$ and translate
into estimates on the difference of the gaps of $D$ and $\Gamma$.
More precisely, we prove the following theorem. Let $\overline{A}$
be the closure of $A$, let $D$ be as in Theorem~\ref{th:1}, and
let $\Theta_{a,b}=\Theta=\{(x,y)\mid 0<x<a,-\frac{b}{a}x+b\leq
y<b\}$, and note that if $D$ is as in Theorem~\ref{th:1} and
contains the points $(a,0)$ and $(0,b)$ then any point in the
first quadrant which is in $(-a,a)\times(-b,b)$ but not in $D$
must be in $\Theta$.
\begin{theorem}\label{th:1.3}
Let the domain $D\subsetneq (-a,a)\times (-b,b)$ be convex and
symmetric about both the $x$-  and $y$-axes. There is a computable
positive function $g_{a,b}=g$ on $\Theta$ such that if $D$
contains $(a,0)$ and $(0,b)$ but $D$ does not contain
$(u_0,v_0)\in \Theta$ then
\begin{equation*}
\lambda_2^{D}-\lambda_1^D \geq
\frac{3\pi^2}{4\max(a,b)^2}+g(u_0,v_0).
\end{equation*}
\end{theorem}
We do not find the largest possible $g$, nor do we know how to do
this. By a computable function we mean  a function of $a$, $b$,
$u_0$, and $v_0$,  involving only elementary one-dimensional
functions. The statement that results if ``computable'' is removed
from the statement of Theorem~\ref{th:1.3} follows fairly quickly
from Theorem~\ref{th:1}. See the end of Section 3.

In Section~4 we discuss possible extensions of
Proposition~\ref{pr:2} in which convexity and symmetry around both
the $x$-  and $y$-axes are respectively replaced by convexity in
$x$ (i.e., two points in $D$ with the same $y$ value can be joined
with a line segment lying in $D$) and symmetry about only the
$y$-axis. We also discuss possible extensions of our results to
certain Schr\"{o}dinger operators and to higher dimensions. These
extensions would lead to inequalities for the difference between
two first eigenvalues but not to inequalities for spectral gaps,
since no analogs of Payne's theorem are known in these settings.

The first use of ratios involving heat kernels to bound gaps was
in \cite{davis}. In \cite{banmen} different proofs of the results
of \cite{davis}  were given as well as a number of interesting
generalizations. (See Section 4 of this paper.) One of these was a
ratio inequality involving integrals of heat kernels. These are
easier to prove than pointwise inequalities for heat kernels and
yield the same information about gaps, which is why we use ratios
of integrals in (\ref{eq:1.1}).

It is easy to modify an example in \cite{smits} to show that if
$H_{\varepsilon} = (-1,1)\times (-1,1)\setminus \{(0,y) \mid
|y|\geq \varepsilon\}$, then the gap of $H_\varepsilon$ goes to
$0$ as $\varepsilon$ approaches $0$. Thus without the convexity
condition or something to replace it, the conclusion of
Theorem~\ref{th:1} does not hold. This example also shows that
Theorems 3.2, 3.3, and 3.4 of \cite{draghici} are incorrect.

\section{Proof of Proposition~\ref{pr:2}}

The proof of Proposition~\ref{pr:2} is based on the connection
between the heat kernel and Brownian motion and the approximation
of Brownian motion by random walks.

In this section we work only with bounded planar domains. Some of
our formulas will hold for all such domains; for these we use
$\Omega$ to designate a domain. Other formulas are not claimed to
hold for all bounded domains but do hold for all bounded convex
domains; for these we use $D$ to designate a domain. We work only
with first, ground state eigenfunctions of a domain $\Omega$, and
we use $\phi^{\Omega}$ to denote this eigenfunction normalized to
integrate to one. The corresponding eigenvalue is denoted by
$\lambda^{\Omega}$.

Standard one-dimensional Brownian motion is denoted by $W_t$,
$t\geq 0$. We use subscripts to denote initial position, as in
$P_x$ and $E_x$, so for example $P_3(W_0=3)=1$. Standard
two-dimensional Brownian motion is denoted $Z_t=(X_t, Y_t)$,
$t\geq 0$. We define $R_0, R_1,\ldots$ to be a random walk such
that $\{R_i-R_{i-1}\}_{i\geq 1}$ are  independent and satisfy
$P(R_i-R_{i-1}=-1)=P(R_i-R_{i-1}=0)=P(R_i-R_{i-1}=1)=1/3$. The
process $W^n$ is the following scaled version of this walk. Let
$R_0=0$. Let $\theta_n=3\cdot 2^{2n-1}$ and note $\var
(2^{-n}R_{\theta_n})=1$. Let $\Theta(n)= \{ k\theta_n^{-1} \mid
k=0,1,2,\ldots\}$. Then $W^n$, started at $0$, is defined by
$W^n_{k\theta_n^{-1}} = 2^{-n}R_k$ so $\var W^n_t =t$, $t\in
\Theta(n)$. For $t\geq 0$ not in $\Theta(n)$, define
$W^n_t=W^n_{k{\theta}^{-n}}$ if
$t\in(k{\theta}_{n}^{-1},(k+1){\theta}_{n}^{-1})$,
$k=0,1,2,\ldots$. Two dimensional scaled random walk is denoted
$Z^n=(X^n,Y^n)$. For this walk, $X^n$ and $Y^n$ are independent
and both have the distribution of $W^n$. We let
$\tau_{\Omega}=\inf\{t>0 \mid Z_t\not\in \Omega\}$ and
$\tau^n_{\Omega}=\inf\{t>0\mid Z^n_t\not\in \Omega\}$, and if $I$
is an interval, $\tau_I=\inf\{t>0 \mid W_t \not\in I\}$ and
$\tau^n_I =\inf \{t>0 \mid W^n_t\not\in I\}$. We use $a_t\sim b_t$
to indicate $\lim_{t\rightarrow \infty} a_t/b_t \in (0,\infty)$.

The eigenfunction expansion of the heat kernel (see Theorem
II.4.13 in \cite{Bass}) implies

\begin{equation}\label{eq:2.1}
P_z\left(\tau_{\Omega} >t\right) = \int_{\Omega}
p_t^{\Omega}(z,w)\,dw \sim e^{-\lambda^{\Omega} t}.
\end{equation}
Theorem~\ref{th:1} follows easily from (\ref{eq:2.1}) and
Proposition~\ref{pr:2}. For a set $A\subset \R^2$ let $A^T=A\cap
\{y>0\}$. Then Payne's theorem implies the minimum of
$\lambda^{D^+}-\lambda^D$ and $\lambda^{D^T}-\lambda^D$ is the gap
of $D$ and the minimum of $\lambda^{\Gamma^+}-\lambda^{\Gamma}$ or
$\lambda^{\Gamma^T}-\lambda^{\Gamma}$ is the gap of $\Gamma$. Now
by (\ref{eq:2.1}),
\begin{equation*}
\lim_{t\to \infty }
e^{(\lambda^{D^+}-\lambda^D)t}\frac{\int_{D^+}p_t^{D^+}(z_0,z)\,dz}{\int_D
p_t^D(z_0,z)\,dz}\in (0,\infty)
\end{equation*}
with a similar formula for $\Gamma^+$ and $\Gamma$. Thus
Proposition~\ref{pr:2} implies
\begin{equation*}
\lambda^{D^+}-\lambda^D >\lambda^{\Gamma^+}-\lambda^{\Gamma}.
\end{equation*}
Rotating $D$ and $\Gamma$  by $90^{\circ}$ and using
Proposition~\ref{pr:2} for these rotated sets gives
\begin{equation*}
\lambda^{D^T}-\lambda^D > \lambda^{\Gamma^T}-\lambda^{\Gamma},
\end{equation*}
and these two inequalities give Theorem~\ref{th:1}.

We note that we can without loss of generality assume that the
closure of the $D$ of Theorem~\ref{th:1} contains $(a,0)$ and
$(0,b)$, since if this is not the case we can replace
$(-a,a)\times(-b,b)$ with the smallest oriented rectangle which
contains $D$, which rectangle could not have a smaller gap. We
also note that the proof of Proposition~\ref{pr:2} is virtually
identical for all $a$ and $b$, and that the assumption that
$(a,0)$ and $(0,b)$ are in $\overline{D}$ does not alter the
proof. Thus we will prove Proposition~\ref{pr:2} only under the
assumptions $\Gamma=S:=(-1,1)\times(-1,1)$ and both $(1,0)$ and
$(0,1)$ belong to $\overline{D}$.

Convex planar domains are Lipschitz, and thus the following lemma,
which states that the heat kernel for convex domains is
intrinsically ultracontractive, is a consequence of the estimates
of \cite{KenigPipher89}. (See the tenth line from the bottom of
page 618 of \cite{KenigPipher89}.) An essentially stronger result
may be found in \cite{BanuelosDavis92}. Intrinsic
ultracontractivity was introduced by Davies and Simon in
\cite{DaviesSimon}.

\begin{lemma}\label{lem:2.1}
There is a positive increasing function $c^D(t)$ on $(0,\infty)$,
and a decreasing function $C^D(t)$ on $(0,\infty)$, such that
$\lim_{t\to \infty} c^D(t)=\lim_{t\to \infty} C^D(t)= 1$,
 and for all $x,y\in
D$ and $t>0$,
\begin{equation}\label{eq:2.2}
\frac{c^D(t)\phi^D(x)\phi^D(y) e^{-\lambda^D t}}{\int_D
\phi(x)^2\,dx} \leq p_t^D(x,y)\leq
\frac{C^D(t)\phi^D(x)\phi^D(y)e^{-\lambda^D t}}{\int_D
\phi(x)^2\,dx}.
\end{equation}
\end{lemma}

Integrating   (\ref{eq:2.2}) in $y$ and using (\ref{eq:2.1}) gives
\begin{equation}\label{eq:4}
\frac{c^D(t)\phi^D(x)e^{-\lambda^D t}}{\int_D \phi(x)^2\,dx} \leq
P_x\left( \tau_D
>t\right) \leq  \frac{C^D(t) \phi^D(x) e^{-\lambda^D t}}{\int_D \phi(x)^2\,dx}.
\end{equation}

One dimensional versions of (\ref{eq:2.2}) and (\ref{eq:4}) where
$D$ is a finite open interval, follow from the references provided
for (\ref{eq:2.2}). These one dimensional inequalities are also
pretty easy to prove directly from the equations of either the
heat kernel or eigenfunctions, which are known for intervals.

The next lemma is a consequence of the classical  fact that given
a two-dimensional Brownian motion $Z_t$, $t\geq 0$, there is a
sequence of $n$-scaled random walks $Z^n_t$, such that for any
$\varepsilon
>0$ and any $K>0$,
\begin{equation}\label{eq:5.1}
P\left(\max_{0\leq s\leq K} |Z_s-Z_s^n| >\varepsilon \right)
\rightarrow 0 \quad\mbox{ as $n\rightarrow \infty$.}
\end{equation}
 A proof of this well known fact is sketched in
\cite{DH}. Here we are abusing notation a little as $Z^n$ in
(\ref{eq:5.1}) stands for a specific $n$-scaled random walk while
in the following lemma it stands for a generic $n$-scaled random
walk.

\begin{lemma}\label{lem:2.2}
Let $\Omega$ be a domain in $\R^2$ and let $Q_1,Q_2,\ldots,Q_m$ be
convex subdomains of $\Omega$. Then for any $0<t_1<\cdots<t_m$,
and any $z\in \Omega$,
\begin{equation}\label{eq:6}
\lim_{n\rightarrow \infty} P_z\left( Z^n_{t_i}\in Q_i, 1\leq i\leq
m, \tau^n_{\Omega} >t_m\right) = P_z\left(Z_{t_i}\in Q_i, 1\leq
i\leq m, \tau_{\Omega}>t_m\right).
\end{equation}
Especially,
\begin{equation}\label{eq:7}
\lim_{n\rightarrow \infty} P_z\left(\tau^n_{\Omega} >t\right) =
P_z\left(\tau_{\Omega} >t\right)\qquad t>0.
\end{equation}
\end{lemma}
 Again, (\ref{eq:6}) and (\ref{eq:7}) are known.
A proof of (\ref{eq:7}) is sketched in \cite{DH}. The equality
(\ref{eq:6}) follows quickly from (\ref{eq:5.1}), the fact that
the probability that $Z_{t_i}$ belongs to the boundary of $Q_i$
equals zero, the fact that the probability that $Z_{t_m}$ belongs
to the boundary of $\Omega$ equals $0$, and the fact (see
\cite{DH}) that the probability that $Z$ hits the boundary of
$\Omega$ for some $t<t_m$ but does not hit the complement of the
closure of $\Omega$ for some $t<t_m$ equals zero.

We denote by $L_t^{\Omega}$, $t\geq 0$, the Markov process which
has transition probabilities
\begin{align}\label{eq:7.1}
l^{\Omega}_t (x,y) &=\frac{p_t^{\Omega}(x,y)
\phi^{\Omega}(y)}{\int_{\Omega}p^{\Omega}_t(x,z)\phi^{\Omega}(z)\,dz}\notag\\
&=p_t^{\Omega}(x,y)\frac{\phi^{\Omega}(y)}{\phi^{\Omega}(x)}e^{\lambda^{\Omega}t},
\end{align}
and stationary distribution
\begin{equation*}
\psi^{\Omega}(x)=\frac{\phi^{\Omega}(x)^2}{\int_\Omega\phi^{\Omega}(y)^2\,dy}.
\end{equation*}
(See Theorem~II.4.13 in \cite{Bass} for (\ref{eq:7.1}).) This
process is often called Brownian motion in $\Omega$ conditioned to
never exit~$\Omega$.

Let $s_1<r<s_2$. The conditional distribution of $Z_r$, $s_1< r<
s_2$, given $Z_{s_1}=z$ and $Z_{s_2}=w$ and $\tau_{\Omega}
>s_2$, is exactly the same as the conditional distribution of
$L_r^{\Omega}$, $s_1< r< s_2$, given $L_{s_1}^{\Omega} =z$ and
$L_{s_2}^{\Omega}=w$. This follows by computing directly the joint
densities of these two processes at points $r_1<r_2<\cdots<r_n$,
where $s_1<r_1$ and $s_2>r_n$, a computation we omit despite the
fact that it is  pleasing to see the eigenfunctions in the
equations for $L$ cancel away.

Now let $s=s(t)$ be such that  $0<s<t$ and both $s$ and $t-2s$
approach infinity as $t$ approaches infinity. Let $\alpha_x(z,w)$
be the joint density of $(Z_s, Z_{t-s})$ given $\tau_D>t$ and
$Z_0=x$. Let $\beta(z,w)$ be the joint density of
$(L^D_s,L^D_{t-s})$ given that $L^D_0$ has density $\psi^D$, its
stationary density. So $\beta$ depends on $t-s$ while $\alpha_x$
depends on $x$, $s$, and $t-s$. Note that $\alpha_x(z,w)$ is the
normalization, to integrate to 1, of
\begin{equation*}
p_s^D(x,z)p_{t-2s}^D(z,w)P_w\left(\tau_D>s\right),
\end{equation*}
while
\begin{equation*}
\beta(z,w) =\psi^D(z)l_{t-2s}^D(z,w).
\end{equation*}
Therefore, (\ref{eq:2.2}) and  (\ref{eq:4}) imply that
$\alpha_x(z,w)/\beta(z,w)$ converges uniformly to 1 in $D\times
D$, as $t\rightarrow \infty$, at a rate which may be taken
independent of $x$.

Thus, if $\mathcal{B}[s,t-s]$ is the $\sigma$-field of the
continuous functions on $[s,t-s]$ generated by the projection
maps, and $x\in D^+$,
\begin{equation}\label{eq:8}
\lim_{t\rightarrow \infty}\sup_{\begin{subarray}{l} A\in
\mathcal{B}[s,t-s]\\ x\in D^+\end{subarray}} \left|
P_x\left(Z_{\cdot}|_{[s,t-s]}  \in A \mid \tau_{D^+}>t\right)-
P_{\psi^{D^+}}\left(L_{\cdot}^{D^+}|_{[s,t-s]}\in
A\right)\right|=0.
\end{equation}

Now let $m$ be an integer which will soon approach infinity, and
let $m'=[\sqrt{m}]$, where $[\quad]$ is the greatest integer
function. Let $(u_0,v_0) \in  \{(x,y)\mid 0<x<1, -x+1\leq y<
1\}\setminus D^{+}$. Put $\Delta_1=D^+\cap \{y<v_0\}$ and
$\Delta_2=D^+\cap \{y>v_0\}$. Let $G(\Delta_1,\Delta_2)$ be the
subset of the continuous functions from $[0,1]$ to $\R^2$ defined
by $G(\Delta_1,\Delta_2)=\{f(0)\in \Delta_1, f(1)\in \Delta_2,
f(t)\in D^+, 0\leq t\leq 1\}$.

Consider the events $F_k=\{ L^{D^+}_{k+\cdot}|_{[0,1]} \in
G(\Delta_1,\Delta_2)\}$, $k=0,1,2,\ldots$. If $L^{D^+}_0$ has
density $\psi^{D^+}$, then $L^{D^+}$ is a stationary process, and
the sequence $I_{F_1},I_{F_2},\ldots,$ is stationary. It is easily
checked that this sequence is ergodic, using the transition
probabilities of $L^{D^+}$ and (\ref{eq:2.2}).  Let
$C_0=P_{\psi^{D^+}}(F_0)$. The ergodic theorem (see
\cite{Breiman}),  which says $\lim_{n\to \infty} \sum_{k=0}^{n-1}
I_{F_k}/n =C_0 \mbox{ a.e.}$, gives--- here we could replace $4/5$
with any constant less than 1---
\begin{equation}\label{eq:9}
\lim_{n\rightarrow \infty}
P_{\psi^{D^+}}\left(\frac{\sum_{k=0}^{n-1} I_{F_k}}{n} >
\frac{4}{5} C_0\right)=1,
\end{equation}
which implies
\begin{equation}\label{eq:10}
\lim_{m\rightarrow\infty}P_{\psi^{D^+}}\left(\frac{\sum_{k=m'}^{m-m'-1}I_{F_k}}{m}>\frac{7}{10}C_0\right)=1.
\end{equation}

Using (\ref{eq:8}), with $m$ and $m'$ in the roles of $t$ and $s$,
(\ref{eq:10}) gives, for $x\in D^+$,
\begin{equation}\label{eq:11}
\lim_{m\rightarrow \infty} P_x\left(\frac{\sum_{k=m'}^{m-m'-1}
I_{\left\{ Z_{k+\cdot}|_{ [0,1]} \in
G(\Delta_1,\Delta_2)\right\}}}{m}
> \frac{3}{5} C_0\mid \tau_{D^+}>m\right)=1.
\end{equation}

We now finish the proof of Proposition~\ref{pr:2}. Recall we are
proving Proposition~\ref{pr:2} only when $\Gamma=S$ and both
$(1,0)$ and $(0,1)$ belong to $\overline{D}$.  Let $(u_0,v_0)\in
S^+\setminus D^+$ and satisfy $(u_0,v_0)\in 2^{-q}\Z^2$, for some
$q=q(u_0,v_0)\in \N$. This is possible since $\cup_n 2^{-n}\Z^2$
is dense in $\R^2$. (This guarantees the discrete walk, for large
enough $n$, hits the line $y=v_0$ when it crosses from below to
above this line.)  Let $z_0=(x_0,y_0)\in D^+$.  For integers
$n\geq q(u_0,v_0)$, and $m$, let $A(n,m)$ be the set of all
sequences $\boldsymbol{y}=(y_0,y_1,\cdots,y_{\theta_n m})$,
satisfying
\begin{equation*}
P_{z_0}\left(Y^n_{k\theta_n^{-1}}=y_k, 0\leq k\leq \theta_n m,
\tau_{D^+}^n>m\right)>0.
\end{equation*}
So certainly  $|y_i-y_{i-1}|$ equals either $2^{-n}$ or zero and
$-1<y_i<1$ for each $i$.

Let $\varepsilon >0$ and let $Q$ be a positive integer. For $m$ a
positive integer let $N(n,Q,m,\varepsilon)$ be the subset of
$A(n,m)$ consisting of those $\boldsymbol{y}$ such that there
exist at least $\varepsilon m$ integers
$h_1<h_2<\cdots<h_{\alpha}$
---so $\alpha\geq \varepsilon m$--- such that $h_1>Q\theta_n$,
$h_{\alpha}<(m-Q)\theta_n$, $|h_i-h_{i-1}|> Q\theta_n$, $2\leq
i\leq \alpha$, and $y_{h_i}=v_0$. Let $m_0=m_0(Q)$ be the smallest
integer such that the probability in (\ref{eq:11}) exceeds $1/2$
if $m\geq m_0$. If $\Y^n_m=(Y^n_{k\theta^{-1}_n})_{ 1\leq k\leq
m\theta_n}$,  we claim that for  $m\geq m_0$ and  $m'>Q$, there
exists $n_1(m)$ such that for $n\geq n_1(m)$,
\begin{equation}\label{eq:12}
P_{z_0}\left( \Y^n_m\in
N\left(n,Q,m,\frac{3C_0}{5(Q+2)}\right)\mid \tau_{D^+}^n>m\right)
>\frac{1}{2}.
\end{equation}
To see this, note that
\[ \left\{\sum_{k=m'}^{m-m'-1}I_{\{Z_k+\cdot|_{[0,1]}\in
G(\Delta_1,\Delta_2)\}}>\frac{3C_0m}{5},\tau_{D^+}>m\right\}\] is
a  union of some of the $2^{m-2m'}$ disjoint sets of the form
\begin{equation*}
  L_j=\left\{Z_i\in
\Delta_{s(i,j)}, m'\leq i\leq m-m'-1,Z_m\in
D^{+},\tau_{D^+}>m\right\}
\end{equation*}
 where $s(i,j)$ is either 1 or 2. Since each $L_j$ is an event of
 the form covered by Lemma~\ref{lem:2.2}, the definition of $m_0$
 and Lemma~\ref{lem:2.2} show that there exits an $n_1(m)$
 such that for $n\geq n_1(m)$,
\begin{equation}\label{eq:12.1}
 P_x\left(\frac{\sum_{k=m'}^{m-m'-1} I_{\left\{ Z^n_{k+\cdot}|_{
[0,1]} \in G(\Delta_1,\Delta_2)\right\}}}{m}
> \frac{3}{5} C_0\mid \tau^n_{D^+}>m\right)>\frac{1}{2}.
\end{equation}

 Now if a
path of $Z^n$ takes values in $\Delta_1$ at $i$ and in $\Delta_2$
at  time $i+1$ then at some time $l(i)\theta_n^{-1}$ between these
two times the $y$-coordinate of the path of $Z^n$ must equal
$v_0$. Furthermore of a collection of $\eta$ such times, each
corresponding to a different integer $i$, at least $\eta/(Q+2)$
may be chosen which are all a distance $Q$ from all the others:
let the smallest be the first chosen, the $Q+2$ smallest be the
second chosen, the $2(Q+2)$ smallest be the third chosen, and so
on.
  Therefore, (\ref{eq:12.1}) and the above observations imply (\ref{eq:12}).

Let $P^{\y}$ designate conditional probability associated with
$Z_0^n,\ldots,Z_m^n$ given $Y_{i\theta_n^{-1}}^n=y_i$, $1\leq i
\leq m\theta_n$.

\begin{lemma}\label{lem:6}
There is an integer $K_0=K_0(u_0)>0$ and a number $0<d=d(u_0)<1$
such that if $\y \in A(n,m)$ and if there are $\lambda$ entries
$y_{j_1},y_{j_2},\ldots,y_{j_\lambda}$ of $\y$ which satisfy
$|{j_i}-{j_{i-1}}|>K_0\theta_n$, $K_0\theta_n <{j_1}$, and
${j_{\lambda}} < (m-K_0)\theta_n$ and $y_{j_i}=v_0$, $1\leq i\leq
\lambda$, then there is an integer $n_2=n_2(m)$ depending only on
$m$ but not on $\y$ such that

\begin{equation}\label{eq:13}
\frac{P^{\y}_{z_0}\left(\tau_{D^+}^n>m\right)}{P^{\y}_{z_0}\left(\tau_D^n>m\right)}
\leq
d^\lambda\frac{P_{z_0}\left(\tau_{S^+}^n>m\right)}{P_{z_0}\left(\tau_S^n>m\right)},\quad
n\geq n_2(m).
\end{equation}
\end{lemma}

Before proving Lemma~\ref{lem:6} we show how it implies
Proposition~\ref{pr:2}. Let $\delta_0=3C_0/(5(K_0+2))$. Now
(\ref{eq:13}) and the definition of $N(n,K_0,m,\delta_0)$ yield
\begin{align*}
P^{\y}_{z_0}\left(\tau_{D^+}^n>m\right)&P_{z_0}\left(\Y^n_m=\y\right)\\
\leq
&d^{\delta_0m}\frac{P_{z_0}\left(\tau^n_{S^+}>m\right)}{P_{z_0}\left(\tau_S^n>m\right)}
\times P_{z_0}^{\y}(\tau_D^n>m)P_{z_0}(\Y^n_m=\y),\quad n\geq
n_2(m),
\end{align*}
for $\y \in N(n,K_0,m,\delta_0)$. Summing these over all $\y$ in
$N(n,K_0,m,\delta_0)$ gives

\begin{align}\label{eq:14}
&\frac{P_{z_0}\left(\tau_{D^+}^n> m, \Y^n_m\in
N(n,K_0,m,\delta_0)\right)}{P_{z_0}\left(\tau_D^n>m,\Y^n_m\in
N(n,K_0,m,\delta_0)\right)}\\
&\leq d^{\delta_0 m}
\frac{P_{z_0}\left(\tau_{S^+}^n>m\right)}{P_{z_0}\left(\tau_S^n>m\right)},
\qquad n\geq n_2(m).\notag
\end{align}

Now (\ref{eq:12}) implies that for $m\geq m_0(K_0)$, $m'>K_0$, and
$n\geq n_1(m)$, the numerator in the left side of (\ref{eq:14}) is
at least half of $P_{z_0}(\tau_{D^+}^n>m)$, while of course the
denominator on the left of (\ref{eq:14}) is no larger than
$P_{z_0}(\tau_D^n>m)$. Thus (\ref{eq:14}) gives
\begin{equation}\label{eq:15}
\frac{P_{z_0}\left(\tau_{D^+}^n>m\right)}{P_{z_0}\left(\tau_D^n>m\right)}
\leq 2 d^{\delta_0 m}
\frac{P_{z_0}\left(\tau_{S^+}^n>m\right)}{P_{z_0}\left(\tau_S^n>m\right)},
\end{equation}
if $m\geq m_0$, $m'>K_0$, and $n\geq n_1(m)$.

 Letting $n \to \infty$ for fixed $m$ and using  (\ref{eq:7})  gives Proposition~\ref{pr:2} in
the case $\Gamma = S$ and $t$ an integer and $o(1) \leq
2d^{\delta_0 t}$, which easily implies Proposition~\ref{pr:2} for
all $t$. We note that, by reasoning similar to that of the
paragraph containing (\ref{eq:2.1}), this bound on $o(1)$ implies
$\lambda^{D^+}-\lambda^D-\delta_0(\log d)\geq \lambda^{S^+}
-\lambda^S$.

The proof of Lemma~\ref{lem:6} requires Lemmas \ref{lem:7},
\ref{lem:8}, and \ref{lem:9} below. Lemma~\ref{lem:7} is from
\cite{davis}. See \cite{DH} for an easier proof.

\begin{lemma}\label{lem:7}
Let $m>0$ be a positive integer, and let $f$ and $g$ be integer
valued functions on $\{0,1,2,\ldots,m\}$ such that $2\leq f(k)\leq
g(k)$, $0\leq k\leq m$. Let $R_0,R_1,\ldots$ be random walk as
defined at the beginning of this section and let $i_0$ be an
integer in $(0,f(0))$. Then
\begin{equation}\label{eq:16}
\frac{P_{i_0}\left(0<R_k<f(k),0\leq k\leq
m\right)}{P_{i_0}\left(0<|R_k|<f(k),0\leq k\leq m\right)} \leq
\frac{P_{i_0}\left(0<R_k<g(k),0\leq k \leq
m\right)}{P_{i_0}\left(|R_k|<g(k),0\leq k \leq m\right)}.
\end{equation}
\end{lemma}

The eigenfunctions of the intervals $(0,1)$ and $(-1,1)$ are
$\phi^{(0,1)}(x)=\frac{\pi}{2}\sin\pi x$ and
$\phi^{(-1,1)}(x)=\frac{\pi}{4}\cos(\frac{\pi}{2}x)$. For
$0<\alpha <1$, let \[\beta(\alpha)=\frac{\int_0^\alpha
\phi^{(0,1)}(x)^2\,dx}{\int_0^1\phi^{(0,1)}(x)^2\,dx}\] and
\[\gamma(\alpha)=\frac{\int_{\{|x|<\alpha\}}\phi^{(-1,1)}(x)^2\,dx}{\int_{-1}^1\phi^{(-1,1)}(x)^2\,dx}.\]Then
since $\phi^{(0,1)}$ and $\phi^{(-1,1)}$ have the same shape (or
by direct calculation) we see that $\beta(\alpha)<\gamma(\alpha)$.

\begin{lemma}\label{lem:8}
Let $0<\alpha<1$ and let $\varepsilon>0$, and put
$\beta=\beta(\alpha)$ and $\gamma=\gamma(\alpha)$. There is a
number $K=K(\alpha,\varepsilon)$ such that, for any integer $m$,
if $t_1,t_2,\ldots,t_m$, and $T$ are numbers such that
\begin{equation}\label{eq:17}
K<t_1<t_2<\cdots<t_m<T-K \mbox{  satisfy $|t_i-t_{i-1}|\geq K$,}
\end{equation}
 and $w\in (0,1)$, then
\begin{align}
P_w\left(W_{t_i}<\alpha,1\leq i\leq m \mid \tau_{(0,1)}>T\right)
& \leq (\beta+\varepsilon)^m,\label{eq:18}\\
\intertext{and} P_w\left(|W_{t_i}|<\alpha, 1\leq i\leq m\mid
\tau_{(-1,1)}>T\right) & \geq (\gamma-\varepsilon)^m.\label{eq:19}
\end{align}
\end{lemma}

{\it Proof~} We prove (\ref{eq:18}). The proof of (\ref{eq:19}) is
similar. If $0<s<t$ and $x\in (0,1)$, the joint density
$\eta(y,z)$ of $(W_s,W_{t})$, conditioned on $\tau_{(0,1)}>t$ and
$W_0=x$, is the normalization of
$p_s^{(0,1)}(x,y)p_{t-s}^{(0,1)}(y,z)$. Thus,
 the density $h_{x,z}$ of $W_s$ given
$W_0=x$ and $W_{t}=z$, and $\tau_{(0,1)}>t$ is
\begin{equation*}
h_{x,z}(y)=\frac{p_s^{(0,1)}(x,y)p_{t-s}^{(0,1)}(y,z)}{\int_0^1p_{s}^{(0,1)}(x,y)p_{t-s}^{(0,1)}(y,z)\,dy}.
\end{equation*}
 The one dimensional version of (\ref{eq:2.2})
 and the fact that $C^{(0,1)}(v)$
decreases to $1$ and $c^{(0,1)}(v)$  increases to $1$ as $v$
increases imply that
 if
$v=\min(s,t-s)$, then for all $z$ and $x$ in $(0,1)$
\begin{align}\label{eq:21.1}
&\frac{e^{-\lambda^{(0,1)}t}c^{(0,1)}(v)^2\phi^{(0,1)}(x)\phi^{(0,1)}(y)^2{}{}\phi^{(0,1)}(z)}{\left(\int_0^1\phi^{(0,1)}(x)^2\,dx\right)^2}\\
{}&{}\leq\frac{e^{-\lambda^{(0,1)}s}e^{-\lambda^{(0,1)} (t-s)} c^{(0,1)}(s) c^{(0,1)}(t-s)\phi^{(0,1)}(x)\phi^{(0,1)}(y)^2\phi^{(0,1)}(z)}{\left(\int_0^1\phi^{(0,1)}(x)^2\,dx\right)^2}\notag\\
{}&{}\leq p_s^{(0,1)}(x,y)p_{t-s}^{(0,1)}(y,z)\notag\\
&{}\leq{}\frac{e^{-\lambda^{(0,1)}s}e^{-\lambda^{(0,1)}(t-s)}C^{(0,1)}(s)C^{(0,1)}(t-s)\phi^{(0,1)}(x)\phi^{(0,1)}(y)^2\phi^{(0,1)}(z)}{\left(\int_0^1\phi^{(0,1)}(x)^2\,dx\right)^2}\notag\\
&{}\leq{}\frac{e^{-\lambda^{(0,1)}t}C^{(0,1)}(v)^2\phi^{(0,1)}(x)\phi^{(0,1)}(y)^2\phi^{(0,1)}(z)}{\left(\int_0^1\phi^{(0,1)}(x)^2\,dx\right)^2}.\notag
\end{align}

Thus
\begin{align}\label{eq:21.2}
\frac{\phi^{(0,1)}(x)\phi^{(0,1)}(z)e^{-\lambda^{(0,1)} t}
c^{(0,1)}(v)^2}{\int_0^1\phi^{(0,1)}(x)^2\,dx} &\leq \int_0^1
p_s^{(0,1)}(x,y)p_{t-s}^{(0,1)}(y,z)\,dy\\
&\leq\frac{\phi^{(0,1)}(x)\phi^{(0,1)}(z)e^{-\lambda^{(0,1)} t}
C^{(0,1)}(v)^2}{\int_0^1 \phi^{(0,1)}(x)^2\,dx}.\notag
\end{align}

Together, (\ref{eq:21.1}) and (\ref{eq:21.2}) imply that
\begin{equation*}
\frac{c^{(0,1)}(v)^2}{C^{(0,1)}(v)^2}\psi^{(0,1)}(y)<h_{x,z}(y)<\frac{C^{(0,1)}(v)^2}{c^{(0,1)}(v)^2}\psi^{(0,1)}(y).
\end{equation*}

This implies that given $\varepsilon >0$ there is
$Q_{(0,1)}=Q_{(0,1)}(\alpha,\varepsilon)$ such that for $v\geq
Q_{(0,1)}/2$ and all $x,z\in (0,1)$,
\begin{equation}\label{eq:20}
\int_0^{\alpha} h_{x,z}(y)\,dy <\beta+\varepsilon.
\end{equation}

Let $\hat{P}$ be conditional probability given $\tau_{(0,1)}>T$.
Then under $\hat{P}$, $W_t$ is a Markov process, although not with
stationary transition probabilities. Let $A_i=\{W_{t_i}<\alpha\}$.
Let $t^+_m=T$, and for $1\leq i <m$ let $t_i^+=(t_i+t_{i+1})/2$.
Suppose $T$ and $t_i$, $1\leq i \leq m$, satisfy (\ref{eq:17})
with $Q_{(0,1)}$ in place of $K$. Then (\ref{eq:18}) holds by the
following  argument, which can be made rigorous by changing
$W_{t_i}=z_i$ to $W_{t_i}\in [z_i,z_i+dz]$.
\begin{align*}
\hat{P}_w\left(A_1\mid W_{t_1^+}=z_1\right) &=
\frac{\hat{P}_w\left(A_1,W_{t_1^+}=z_1\right)}{\hat{P}_w\left(W_{t_1^+}=z_1\right)}\\
&=\frac{P_w\left(A_1,W_{t_1^+}=z_1,\tau_{(0,1)}>t_m^+\right)}{P_w\left(W_{t_1^+}=z_1,\tau_{(0,1)}>t_m^+\right)}\\
&=
\frac{P_w\left(A_1,W_{t_1^+}=z_1,\tau_{(0,1)}>t_1^+\right)P_{z_1}\left(\tau_{(0,1)}>t_m^+-t_1^+\right)}{P_w\left(W_{t_1^+}=z_1,\tau_{(0,1)}>t_1^+\right)P_{z_1}\left(\tau_{(0,1)}>t_m^+-t_1^+\right)}\\
&<\beta+\varepsilon \qquad\mbox{ by (\ref{eq:20}).}
\end{align*}
 Using similar ratios and the Markov property we obtain $\p_w(A_2\mid
A_1,W_{t_1^+}=z_1,W_{t_2^+}=z_2)=\p_{z_1}(A_1\mid
W_{t_1^+}=z_2)<\beta+\varepsilon$ so that $\p_w(A_1\cap A_2 \mid
W_{t_2^+}=z_m)< (\beta+\varepsilon)^2$. Proceeding in this manner
gives
\begin{equation*}
\p_w\left(A_1\cap A_2\cap \cdots\cap A_m \mid W_{t_m^+}=z_m\right)
< (\beta+\varepsilon)^m,
\end{equation*}
and integrating over $z_m$ gives (\ref{eq:18}). Similarly there is
a $Q_{(-1,1)}(\alpha,\varepsilon)$ such that if $T$ and $t_i$,
$1\leq i\leq m$, satisfy (\ref{eq:17}) with $Q_{(-1,1)}$ in place
of $K$, then (\ref{eq:19}) holds. So $K(\alpha,\varepsilon)$ can
be and is taken to be the smallest integer larger than
$\max(Q_{(0,1)}(\alpha,\varepsilon),Q_{(-1,1)}(\alpha,\varepsilon))$.

Let $d=d(u_0)=(\frac{\beta(u_0)}{\gamma(u_0)}+1)/2<1$, and let
$\varepsilon=\varepsilon(u_0)$ satisfy
$\frac{\beta+\varepsilon}{\gamma-\varepsilon}<
\frac{1}{2}({d+\frac{\beta(u_0)}{\gamma(u_0)}})$. Let
$K_0=K(u_0,\varepsilon(u_0))$. Lemma~\ref{lem:8} implies
\begin{align}\label{eq:20.1}
&\frac{P_w\left(W_{t_i}<u_0, 1\leq i\leq m,\mbox{ and }
\tau_{(0,1)}>T\right)}{P_w\left(|W_{t_i}|<u_0, 1\leq i \leq m,
\mbox{ and } \tau_{(0,1)}>T\right)} \\&< d^m
\frac{P_w\left(\tau_{(0,1)}>T\right)}{P_w\left(\tau_{(-1,1)}>T\right)}.\notag
\end{align}

\begin{lemma}\label{lem:9}
Let  $T$, $t_1,t_2,\ldots,t_m$ be as in Lemma~\ref{lem:8}, and
suppose in addition that all $t_i$, $1\leq i\leq m$, are in
$\Theta(n)$ for some $n$ and $w=l2^{-n}$ where $l$ is an integer
such that $0<l<2^n$. Then there is an integer $N(m,T)$ such that
\begin{align}\label{eq:21}
&\frac{P_w\left(W_{t_i}^{(n)}<u_0, 1\leq i\leq m,\mbox{ and }
\tau_{(0,1)}^{(n)}>T\right)}{P_w\left(|W_{t_i}^{(n)}|<u_0, 1\leq i
\leq m, \mbox{ and } \tau_{(0,1)}^{(n)}>T\right)} \\&< d^m
\frac{P_w\left(\tau_{(0,1)}^{(n)}>T\right)}{P_w\left(\tau_{(-1,1)}^{(n)}>T\right)}\quad
n\geq N(m,T).\notag
\end{align}
\end{lemma}
{\it Proof} That the limit as $n\rightarrow \infty$ of both of the
numerators and both of the denominators in (\ref{eq:21}) exists
follows from a one dimensional version of Lemma~\ref{lem:2.2},
which implies that these limits are the analogous probabilities
for Brownian motion $W_t$. This one dimensional version follows
from the classical result of Skorohod that processes with the
distribution of $W^{(n)}$ may be embedded in $W$ in such a way
that given $t>0$ and $\varepsilon>0$,
$P(|W_s-W_s^{(n)}|<\varepsilon, 0\leq s\leq t)>1-\varepsilon$ for
large enough $n$. Together with (\ref{eq:20.1}) this establishes
Lemma~\ref{lem:9}.

We note that, for $z_0=(x_0,y_0)$, the independence of the
components $X^n$ and $Y^n$ of $Z^n$ implies that if $m$ is an
integer both
\begin{align*}
P_{z_0}\left(\tau_{S^+}^n>m\right) =&P_{x_0}\left(\inf\left\{t\mid
W_t^n\not\in
(0,1)\right\}>m\right)\\
\times &P_{y_0}\left(\inf\left\{t\mid
W_t^n\not\in(-1,1)\right\}>m\right)
\end{align*}
and
\begin{align*}
  P_{z_0}\left(\tau_S^n>m\right)=&
P_{x_0}\left(\inf\left\{t\mid W_t^n\not\in(-1,1)\right\}>m\right)\\
&\times P_{y_0} \left(\inf\left\{t\mid
W_t^n\not\in(-1,1)\right\}>m\right)
\end{align*}
where $W_t^n$ is the one dimensional random walk  defined at the
beginning of this section. Thus the ratio of the right hand side
of (\ref{eq:13}) is
\begin{equation}\label{eq:22}
\frac{P_{2^nx_0}\left(R_k\in(0,2^n),0\leq k \leq
m\theta_n\right)}{P_{2^nx_0}\left(R_k\in(-2^n,2^n),0\leq k\leq
m\theta_n\right)}=\frac{P_{x_0}\left(\tau_{(0,1)}^n>m\right)}{P_{x_0}\left(\tau_{(-1,1)}^n>m\right)}.
\end{equation}

Now we finish the proof of Lemma~\ref{lem:6}. We think of $\y$ as
fixed. Let $q_n(k)=\max\{i\in \Z\mid(i2^{-n},y_k)\in D\}$, $1\leq
k\leq m2^n$, and $Q(n)=\max\{i\in \Z\mid(i2^{-n},v_0)\in D\}$, and
put $\hat{q}_n(k)=Q(n)$ for $k= {j_i}$, $1\leq i\leq \lambda$,
$\hat{q}_n(k)=2^n$ for all other $k\in \N$. Note $q_n(k)\leq
\hat{q}_n(k)$, $0\leq k\leq m\theta_n$.

Then using Lemma~\ref{lem:7} we get
\begin{equation}\label{eq:23}
\frac{P_{2^n x_0}\left(0<R_k<q(k), 1\leq k\leq
m\theta_n\right)}{P_{2^n x_0}\left(|R_k|<q(k), 1\leq k\leq
m\theta_n\right)}\leq \frac{P_{2^n x_0}\left(0<R_k<\hat{q}(k),
1\leq k\leq m\theta_n\right)}{P_{2^n x_0}\left(|R_k|<\hat{q}(k),
1\leq k\leq m\theta_n\right)}.
\end{equation}

Now, upon scaling, the left hand ratio in (\ref{eq:23}) becomes
\begin{equation*}
\frac{P^{\y}_{x_0}(\tau_{D^+}^n>m)}{P^{\y}_{x_0}(\tau_D^n>m)},
\end{equation*}
which is the left hand side of (\ref{eq:13}), while the ratio on
the right hand side of (\ref{eq:23}) becomes
\begin{equation*}
\frac{P_{x_0}\left(0<W^n_{{j_i}\theta_n^{-1}}<u_0,1\leq
i\leq\lambda,\tau_{(0,1)}^n>m\right)}{P_{x_0}\left(|W^n_{{j_i}\theta_n^{-1}}|<u_0,
1\leq i \leq \lambda, \tau_{(-1,1)}^n>m\right)},
\end{equation*}
which by Lemma~\ref{lem:9} does not exceed
\begin{equation*}d^\lambda
\frac{P_{x_0}(\tau_{(0,1)}^n>m)}{P_{x_0}(\tau_{(-1,1)}^n>m)},
\end{equation*}
which is the right hand of (\ref{eq:13}). This proves
(\ref{eq:13}).

\section{Proof of Theorem~\ref{th:1.3}}
Again we just prove Theorem~\ref{th:1.3} in the case $\Gamma=S$
and when $(0,1)$ and $(1,0)$ are in $\overline{D}$ and we let
$(u_0,v_0)$ be a point in $(S\setminus D)\cap \{x>0,y>0\}$. Again
the argument for arbitrary $a$ and $b$ is virtually identical.

In Section~2, just after (\ref{eq:15}), we showed that
\begin{equation}\label{eq:24}
\lambda^{D^+}-\lambda^D-\delta_0 \log d \geq
\lambda^{S^+}-\lambda^S.
\end{equation}
Now $d=d(u_0)$ is explicitly defined in the last section.
 To complete the proof of Theorem~\ref{th:1.3}
 we show that there is a computable positive number
$h(u_0,v_0)$ depending only on $u_0$ and $v_0$ alone for which the
inequality that results when that number is substituted for
$\delta_0$ in (\ref{eq:24}) is true. In other words,
\begin{equation}\label{eq:25}
\lambda^{D^+}-\lambda^D-h(u_0,v_0)\log d(u_0)\geq
\lambda^{S^+}-\lambda^S,
\end{equation}
and Theorem~\ref{th:1.3} in the case $\Gamma=S$ is verified, where
\[g(u,v)=\min(h(u,v)\log d(u), h(v,u)\log d(v)).\] Here we use the
fact that if $(u_0,v_0)$ is omitted from $D^+$ then $(v_0,u_0)$ is
omitted from the $90^\circ$ rotation of $D^T$.

We need to show how to produce a positive function of $u_0$ and
$v_0$ which is smaller than $\delta_0$.
  Recall
$\delta_0=3P_{\psi^{D^+}}(F_0)/(5(K_0+2))$. Now the proof of
Lemma~\ref{lem:10} will show how to bound $P_{\psi^{D^+}}(F_0)$
below.
 Also, $K_0=K(u_0,\varepsilon(u_0))$, and $\varepsilon(u_0)$ is easily
computed while $K_0$ is defined just above (\ref{eq:20.1}). The
proof of Lemma~\ref{lem:8} shows that if we can produce explicit
versions of the functions $c^{(0,1)}(v)$ and $C^{(0,1)}(v)$ then
the desired explicit upper bound for $K_0$ can be achieved.
 These explicit versions can be found using either the exact formulas for the
 one-dimensional heat kernels or the exact formulas of the eigenvalues
of an interval.

\begin{lemma}\label{lem:10} The quantity $P_{\psi^{D^+}}(F_0)$ may be bounded below by
a positive number which depends only on $u_0$ and $v_0$.
\end{lemma}

We first prove some estimates for $p_1^{D^+}(z,w)$ and
$\phi^{D^+}$. For $z\in D^+$, let $\eta(z)$ be the distance from
$z$ to the boundary of $D^+$, and let $D_{\varepsilon}^+$ be all
points $z$ of $D^+$ which satisfy $\eta(z)>\varepsilon$ . For $z$,
$w\in D^+$, let $R(z,w)$ be the rectangle which lies in $D^+$, has
one side of length $\min(\eta(z),\eta(w))$,  and contains $z$ and
$w$ and such that both $z$ and $w$ are  a distance
$\min(\eta(z),\eta(w))/2$ from three of the four sides of
$R(z,w)$. (So, if either $z$ or $w$ are close to the boundary of
$D^+$ and $z$ and $w$ are far apart, $R(z,w)$ is a long skinny
rectangle and $z$ and $w$ are both close to short sides of
$R(z,w)$.) The convexity of $D^+$ guarantees $R(z,w)\subset D^+$.
The heat kernel of a rectangle is exactly known and using this
exact formula and the fact that the diameter of $D^+$ does not
exceed 3, so the long side of $R(z,w)$ is no longer than 3, it is
easy to give the equation of a positive increasing function
$\beta(s)$, $s>0$ such that
\begin{equation*}
p_1^{R(z,w)}(z,w)\geq \beta(\min(\eta(z),\eta(w))), \quad z,w\in
D^{+},
\end{equation*}
which implies
\begin{equation}\label{eq:26}
p_1^{D^+}(z,w)\geq \beta(\min(\eta(z),\eta(w))),\quad z,w\in
D^{+}.
\end{equation}
Also, it is immediate that
\begin{equation}\label{eq:27}
p_1^{D^+}(z,w)\leq p_1^{\R^2}(z,w)=\frac{1}{2\pi}.
\end{equation}

Now since the rectangle $(0,1/2)\times(-1/2,1/2)$, which has first
eigenvalue $3\pi^2$, is contained in $D^+$, we know
\begin{equation}\label{eq:28}
\lambda^{D^+}\leq 3\pi^2.
\end{equation}
Using (\ref{eq:27}) and (\ref{eq:28}) we get
\begin{equation*}
e^{-3\pi^2}\phi^{D^+}(z)\leq e^{-\lambda^{D^+}}
\phi^{D^+}(z)=\int_{D^+}\phi^{D^+}(x)p_1^{D^+}(x,z)\,dz\leq
\frac{1}{2\pi}, \quad z\in D^+,
\end{equation*}
yielding
\begin{equation}\label{eq:29}
\phi^{D^+}(z)\leq \frac{e^{3\pi^2}}{2\pi},\quad z\in D^+.
\end{equation}
Also,
\begin{align}\label{eq:30}
\int_{D^+}\phi^{D^+}(z)^2\,dz &= \mbox{ area }D^+
\int_{D^+}\phi^{D^+}(z)^2\,d\left(\frac{z}{\mbox{area
$D^+$}}\right)\\
&\geq \mbox{ area
$D^+$}\left[\int_{D^+}\phi^{D^+}(z)\,d\left(\frac{z}{\mbox{area
$D^+$}}\right)\right]^2\notag\\
&=(\mbox{area $D^+$})\cdot (\mbox{area $D^+$})^{-2} \geq
\frac{1}{2}.\notag
\end{align}
 Now
if $z\in D^+\setminus D^+_\varepsilon$, the convexity of $D^+$
implies that either there is a point on the vertical line through
$z$ which belongs to the boundary of $D^+$ and is a distance at
most $\sqrt{2}\varepsilon$ from $z$ or a point on the horizontal
line through $z$ which belongs to the boundary of $D^+$ and is a
distance at most $\sqrt{2}\varepsilon$ from $z$. Let
$L^x=\{(x,y)\mid -\infty <y<\infty\}$ and $L_y=\{(x,y)\mid -\infty
<x <\infty\}$. Let $\partial D$ stand for the boundary of $D$.
Then for almost every $x$ in $(0,1)$ and every $y$ in $(-1,1)$,
$L^x\cap
\partial D^+$ and $L_y\cap \partial D^+$ each consist of exactly
two points. Let $L^x(a)$ be all points $w$ in $D^+$ which belong
to $L^x$ and such that there is $z\in L^x\setminus D^+$ such that
$|z-w|\leq a$. Let $L_y(a)$ be all points $w$ in $D^+$ that belong
to $L_y$ and such that there is a point $z$ in $L_y\setminus D^+$
such that $|z-w|\leq a$. Then the length of $L_y(a)$ does not
exceed $2a$ when $L_y\cap \partial D^+$ consists of two points,
with a similar inequality for $L^x(a)$. Since
\[ \bigcup_{0\leq x \leq 1} L^x(\sqrt{2}\varepsilon) \cup
\bigcup_{-1\leq y\leq 1} L_y(\sqrt{2}\varepsilon)\supset
D^+\setminus D^+_\varepsilon,\]
\begin{equation}\label{eq:31}
|D^+\setminus D^+_\varepsilon|_2\leq
\int_{0}^{1}|L^x(\sqrt{2}\varepsilon)|_1\,dx+\int_{-1}^1|L_y(\sqrt{2}\varepsilon)|_1\,dy\leq
3\cdot 2\sqrt{2}\varepsilon,
\end{equation}
where $|\quad|$ is respectively Lebesgue two and one dimensional
measure.

The inequalities (\ref{eq:29}), (\ref{eq:30}), and (\ref{eq:31})
give
\begin{equation}\label{eq:32}
\int_{D^+\setminus D^+_\varepsilon}\psi^{D^+}\leq
6\sqrt{2}\varepsilon\left(\frac{e^{3\pi^2}}{2\pi}\right)^2\cdot2.
\end{equation}
Thus, since $\int_{D^+}\psi^{D^+}(z)\,dz=1$, we have that if
$\gamma=[6\sqrt{2}(e^{3\pi^2}/2\pi)^2\cdot 2]^{-1}\cdot1/2$, then
$\int_{D^+\setminus D^+_\gamma} \psi^{D^+} \leq 1/2$, which,
together with the symmetry about the $x$-axis of $\phi^{D^+}$ and
thus $\psi^{D^+}$, and the fact that $\int_{D^+}\psi^{D^+}=1$
gives
\begin{equation}\label{eq:33}
\int_{D_{\gamma}^+\cap \{y<0\}}\psi^{D^+} \geq \frac{1}{4}.
\end{equation}

Furthermore, by (\ref{eq:29}) and (\ref{eq:31}),
\begin{equation}\label{eq:33.1}
\int_{D^+_\gamma}\phi^{D^+}(x)\,dx\geq \frac{1}{4}.
\end{equation}

We also note that, by (\ref{eq:26}) and (\ref{eq:33.1}),
\begin{align}\label{eq:34}
\phi^{D^+}(z) &=
e^{\lambda^{D^+}}\int_{D^+}\phi^{D^+}(x)p_1^{D^+}(x,z)\,dx\\
&\geq \int_{D^+}\phi^{D^+}(x)p_1^{D^+}(x,z)\,dx\notag\\
&\geq
\int_{D^+_\gamma}\phi^{D^+}(x)p_1^{D^+}(x,z)\,dx\notag\\
&\geq\frac{1}{4}\beta(\min(\eta(z),\gamma)).\notag
\end{align}

We now bound $P_{\psi^{D^+}}(F_0)$ from below. Using the
transition probabilities for $L_t^{D^+}$ given by (\ref{eq:7.1}),
we have

\begin{equation*}
P_{\psi^{D^+}}(F_0) = \int_{D^+\cap
\{y<v_0\}}\left[\int_{D^+\cap\{y>v_0\}}l_1^{D^+}(x,y)\,dy\right]\psi^{\Omega}(x)\,dx.
\end{equation*}
Now the open triangle $T(v_0)$ with vertices $(0,1)$, $(0,v_0)$,
and $(1-v_0,v_0)$ must lie in $D^+\cap\{y>v_0\}$. Let $t(v_0)$ be
the middle third of this open triangle, that is, $t(v_0)$ is the
translation of $\frac{1}{3}T(v_0)$ satisfying that the medians of
$t(v_0)$ and $T(v_0)$ meet in the same place. Then all points of
$t(v_0)$ are at least  a distance $v_0/3$ from $\partial {D^+}$,
and so since the area of $t(v_0)$ equals $v_0^2/18$,
\begin{align*}
P_{\psi^{D^+}}(F_0) &\geq \int_{ D^+_{\gamma}\cap \{y<v_0\}}
\left[\int_{t(v_0)}l_1^{D^+}(x,y)\,dy\right]\psi^{D^+}(x)\,dx\\
&\geq \int_{ D^+_{\gamma}\cap
\{y<v_0\}}\left[\int_{t(v_0)}p_1^{D^+}(x,y)\frac{2\pi\beta(\min(\gamma,\frac{v_0}{3}))}{4 e^{3\pi^2}}\,dy\right]\psi^{D^+}(x)\,dx\\
& \geq
\frac{2\pi\beta(\min(\gamma,\frac{v_0}{3}))}{4e^{3\pi^2}}\int_{D^+_{\gamma}\cap
\{y<0\}}\left[\int_{t(v_0)}p_1^{D^+}(x,y)\,dy\right]\psi^{D^+}(x)\,dx\\
& \geq
\frac{2\pi\beta(\min(\gamma,\frac{v_0}{3}))}{4e^{3\pi^2}}\int_{D^+_{\gamma}\cap
\{y<0\}}
\beta\left(\min(\gamma,\frac{v_0}{3})\right)\frac{v_0^2}{18}\psi^{D^+}(x)\,dx\\
& \geq \frac{2\pi\beta(\min(\gamma,\frac{v_0}{3}))^2}{4e^{3\pi^2}}
\cdot\frac{v_0^2}{18}\cdot \frac{1}{4},
\end{align*}
using (\ref{eq:29}) and (\ref{eq:34}) in the second inequality,
(\ref{eq:26}) in the fourth, and (\ref{eq:33}) in the last. This
proves Lemma~\ref{lem:10}, which completes the proof of
Theorem~\ref{th:1.3}.

The existence of a positive function $g(u,v)$ which satisfies the
statement that results if ``computable'' is removed from the
statement of Theorem~\ref{th:1.3} follows from Theorem~\ref{th:1}.
For this existence is equivalent to the statement that the
supremum of the set of all gaps of convex doubly symmetric domains
contained in $\Gamma$ with closures containing $(a,0)$ and $(0,b)$
which do not contain $(u_0, v_0)$ is less than the gap of
$\Gamma$. Suppose by way of contradiction that this is not the
case. Pick a sequence of these domains such that their gaps
converge to the gap of $\Gamma$. Pick a subsequence $D_n$ of these
domains such that $\sup\{y\mid (x,y)\in D_n\}$ converges, say to
$f(x)$, for each rational $x$. Let $D'$ be the unique doubly
symmetric convex domain  such that $\sup\{(x,y)\mid y\in
D'\}=f(x)$ for each rational $x$ in $(-1,1)$. The monotonicity
property of the eigenvalues (see Theorem VI.3 in
\cite{CourantHilbert1}) together with the fact that for
$\varepsilon>0$ if $n$ is large enough then
$(1-\varepsilon)D_n\subset D'\subset (1+\varepsilon)D_n$ implies
that the gap of $D'$ equals the gap of $\Gamma$. And $(u_0, v_0)$
does not belong to $D'$, yielding a contradiction to
Theorem~\ref{th:1}.

\section{Final Comments}
In this section we drop our convention that $D$ always stands for
a convex domain, and discuss some possible extensions of the
results we have proved. We believe that our proof of
Proposition~\ref{pr:2} will extend fairly easily to prove a
stronger version of Proposition~\ref{pr:2} in which the convexity
and double symmetry of $D$ is weakened to convexity in $x$ and
symmetry about the $y$-axis. Under this weakened condition $D$ is
no longer intrinsically ultracontractive. However the analog
of~(\ref{eq:8}), the only place we used two dimensional intrinsic
ultracontractivity, actually holds for all bounded domains.

In \cite{davis} and \cite{bankro} a higher dimensional analog is
proved of the two-dimensional result of \cite{davis} that
$\lambda^{D^+}-\lambda^D\geq \lambda^{\Gamma^+}-\lambda^{\Gamma}$
if $D$ is simply connected, bounded, symmetric about the $y$-axis,
convex in $x$, and contained in $(-a,a)\times(-b,b)$.  We believe
that properly formulated analogs of Proposition~\ref{pr:2} are
true, but note that, as heat kernels do not see a line in higher
dimensions, care needs to be taken in the formulation.

In \cite{banmen} and later in  \cite{bankro}, \cite{DH} and
\cite{draghici}, theorems related to the results of \cite{davis}
were proved for Schr\"{o}dinger operators with potential kernels
which are symmetric about the $y$-axis and nondecreasing in $x$
for positive $x$ for fixed $y$ and which are defined on bounded
domains symmetric about the $y$-axis. We believe that the methods
in \cite{DH} and those of this paper can be used to prove analogs
of Proposition~\ref{pr:2} of this paper for such operators.


\begin{thebibliography}{99}
\bibitem{BanuelosDavis92}
Rodrigo Ba{\~n}uelos   and Burgess Davis: A geometrical
characterization of intrinsic ultracontractivity for planar
domains with boundaries given by the graphs of functions. Indiana
Univ. Math. J.~{\bf 41}~885--913~(1992)
\bibitem{bankro}
Rodrigo Ba{\~n}uelos   and Pawel Kr{\"o}ger: Gradient estimates
for the ground state {S}chr\"odinger eigenfunction and
applications. Comm. Math. Phys.~{\bf 224}~545--550~(2001)
\bibitem{banmen}
 Rodrigo  Ba{\~{n}}uelos and~Pedro~J. M{\'{e}}ndez-Hern{\'{a}}ndez: Sharp inequalities for heat kernels of {S}chr{\"{o}}dinger
operators and applications to spectral gaps.~J. Funct. Anal.~{\bf
176}~368--399~(2000)
\bibitem{Bass}
Richard F. Bass: Probabilistic techniques in analysis.
Springer-Verlag, New York, 1995
\bibitem{vdb}
 M. van den Berg: On condensation in the free-boson gas and the spectrum of the
Laplacian.~J. Statist. Phys.~{\bf 31}~623--637~(1983)
\bibitem{Breiman}
Leo Breiman:~Probability. Society for Industrial and Applied
 Mathematics, Philadelphia, 1992, corrected reprint of the 1968
 original
\bibitem{CourantHilbert1}
R. Courant and D. Hilbert:~ Methods of mathematical physics.
Volume~I. Wiley Classics Edition, 1989
\bibitem{DaviesSimon}
E. B. Davies and B. Simon: Ultracontractivity and the heat kernel
for Schr\"{o}dinger operators and Dirichlet Laplacians.~J. Funct.
Anal.~{\bf 59}~335--395~(1984)
\bibitem{davis}
Burgess  Davis:~On the spectral gap for fixed membranes.~Ark.
Mat.(1)~{\bf 39}~65--74~(2001)
\bibitem{DH}
 Burgess  Davis and Majid Hosseini:~On ratio inequalities for heat content.~J. London. Math.
Soc. (2)~{\bf 69}~97--106~(2004)
\bibitem{draghici}
Cristina Draghici:~Rearrangement inequalities with application to
ratios of heat kernels.~Potential Anal.~{\bf 22}~351--374~(2005)
\bibitem{KenigPipher89}
Carlos E. Kenig  and Jill Pipher:~The {$h$}-path distribution of
the lifetime of conditioned {B}rownian motion for nonsmooth
domains:~Probab. Theory Related Fields~{\bf 82}~615--623~(1989)
\bibitem{Ling93}
Jun Ling:~A lower bound for the gap between the first two
eigenvalues of {S}chr\"odinger operators on convex domains in
{$S\sp n$} or {${\bf R}\sp n$}.~Michigan Math. J. ~{\bf
40}~259--270~(1993)
\bibitem{payne}
Lawrence E. Payne:~On two conjectures in the fixed membrane
eigenvalue problem.~Z. Angew. Math. Phys.~{\bf 24}~721--729~(1973)
\bibitem{swyy}
I. M. Singer, Bun Wong, Shing-Tung Yau, and Stephen S.-T. Yau: An
estimate of the gap of the first two eigenvalues in the
{S}chr\"odinger operator. Ann. Scuola Norm. Sup. Pisa Cl. Sci.
(4)~{\bf 12}~319--333~(1985)
\bibitem{smits}
Robert G. Smits:~Spectral gaps and rates to equilibrium for
diffusions in convex domains.~Michigan Math. J.~{\bf
43}~141--157~(1996)
\bibitem{Wang00}
Feng-Yu Wang: On estimation of the {D}irichlet spectral
gap.~{Arch. Math. (Basel)}~{\bf 75}~450--455~(2000)
\bibitem{YuZhong86}
Qi Huang Yu and Jia Qing Zhong:~Lower bounds of the gap between
the first and second eigenvalues of the {S}chr\"odinger
operator.~Trans. Amer. Math. Soc.~{\bf 294}~341--349~(1986)
\end{thebibliography}
\end{document}